\documentclass{article} 

\usepackage[all]{xy}
\usepackage{amsmath}
\usepackage{amssymb}
\usepackage{amsthm}

\newtheorem{theorem}{Theorem}[section]
\newtheorem{lemma}[theorem]{Lemma}
\newtheorem{corollary}[theorem]{Corollary}
\newtheorem{proposition}[theorem]{Proposition}

\newtheorem{definition}[theorem]{Definition}

\newtheorem{remark}[theorem]{Remark}
\newtheorem{example}[theorem]{Example}

\def\N{{\mathbb N}}
\def\P{{\mathbb P}}
\def\p{{\mathbb P}}

\def\Z{{\mathbb Z}}

\def\cG{{\cal G}}

\def\cO{{\cal O}}

\def\operatorname#1{\mathop{\rm #1}\nolimits}

\def\deg{\operatorname{deg}}
\def\Bs{\operatorname{Bs}}

\def\sga{{\langle}}
\def\sgc{{\rangle}}

\title{Tangential projections and secant defective varieties}
\author{Roberto Mu\~noz, Jos\'e Carlos Sierra and Luis E. Sol\'a Conde}

\begin{document}

\maketitle

\noindent {\bf Abstract:}   Going one step further in Zak's 
classification of Scorza varieties with secant defect equal to one, we
characterize the Veronese embedding of $\P^n$ given by the complete
linear system of quadrics and its smooth projections from a point as
the only smooth irreducible complex and non-degenerate projective
subvarieties of $\P^N$ that can be projected isomorphically into
$\P^{2n}$ when $N\geq\binom{n+2}{2}-2$.

\section{Introduction}\label{sec:intro}

Any smooth complex projective variety $X\subset\P^N$ of dimension $n$
such that $N \geq 2n+1$ can be projected isomorphically into
$\P^{2n+1}$ by simply choosing a center of projection not meeting the
closure of the secant lines to $X \subset \P^N$. As usual in
projective geometry, associated to this general property a problem of
classification appears: to find the complete list of those
non-degenerate smooth complex projective varieties of dimension $n$
that can be projected isomorphically into $\P^{2n}$. This problem is
solved in low dimension. A smooth non-degenerate curve $C\subset\P^N$
($N\geq 3$) cannot be projected isomorphically onto a plane curve. For
$n=2$ a complete list of surfaces with this property was achieved by
Severi in \cite{severi}:

\begin{theorem}\label{Severi}
  Let $X\subset\P^5$ be a smooth irreducible complex and
  non-degenerate projective surface. If $X$ can be projected
  isomorphically into $\P^4$, then $X$ is the Veronese surface
  $v_2(\P^2)\subset\P^5$.
\end{theorem}

The case $n=3$ was first considered by Scorza in
\cite{scorza} and completed by Fujita in \cite{fujita} (see Theorem
\ref{thm:severi-scorza}). When $n=4$ only some partial results are
known, see \cite{scorza2} and \cite{fu-ro}, where an infinite list of
examples is shown. Hence the problem of getting a complete
classification for arbitrary dimension seems far from being reached.
However, if $N$ is big enough, Zak's Theorem on Scorza varieties
\cite[Ch. VI]{Zak} shows that Severi's Theorem can be generalized in
the following way:

\begin{theorem}\label{thm:zak-severi}
  Let $X\subset\P^N$ be a smooth irreducible complex and
  non-degenerate projective variety of dimension $n\geq 2$. Let $N(n)$
  be ${n+2 \choose 2}-1$. If $X$ can be projected isomorphically into
  $\P^{2n}$, then $N\leq N(n)$ with equality if and only if $X$ is the
  second Veronese embedding $v_2(\P^n)\subset\P^{N(n)}$.
\end{theorem}

The main result in the paper is an extension of this theorem,
conjectured in \cite{ASU}, where a similar statement was proved for
subvarieties of grassmannians of lines:

\begin{theorem}\label{thm:USA}  Let $X\subset\P^N$ be a smooth irreducible
  complex and non-degenerate projective variety of dimension $n$ and
  let $N\geq max\{N(n)-1, 2n+1\}$. If $X$ can be projected
  isomorphically into $\P^{2n}$ then one of the following holds:

\begin{itemize}
\item[{\rm (a)}] $X=v_2(\P^n)\subset\P^{N(n)}$;
\item[{\rm (b)}] $X$ is either the isomorphic projection of
  $v_2(\P^{n})$ into $\P^{N(n)-1}$ or its inner projection
  $B^n\subset\P^{N(n)-1}$.
\end{itemize}
\end{theorem}

Let us observe that, as noted above, the center of an isomorphic
projection of $X \subset \P^{2n+1}$ into $\P^{2n}$ cannot intersect
the secant variety of $X$. Hence the property of being projectable
isomorphically into $\P^{2n}$ is equivalent to the fact that the
dimension of the secant variety is smaller than the expected one
(which is $2n+1$). A variety $X \subset \P^N$ with this property on
the dimension of its secant variety is called {\it $1$-defective} and
the difference between the actual dimension and the expected one is
called the {\it $1$-secant defect of} $X \subset \P^N$. By Terracini's
Lemma (cf. Lemma \ref{lem:terracini}) the dimension of the secant
variety can be computed by looking at the linear space spanned by two
general projective tangent spaces to $X$. This shows that
$1$-defectivity corresponds to the fact that the linear projection of
$X$ from a general tangent space (the {\it tangential projection}) is
not of maximal rank. It is now when tangential projections enter into
the picture. Concretely, Theorem \ref{thm:USA} is an immediate
consequence of the following result:

\begin{theorem}\label{thm:main}
  Let $X\subset\P^N$ be a smooth irreducible complex and
  non-degenerate projective variety of dimension $n$. Assume $X$ is
  $1$-defective and consider $k<k_0$ a positive integer, where $k_0$
  is the least integer verifying $S^{k_0}X=\P^N$.
\begin{enumerate}
\item[{\rm (a)}] If the general $k$-tangential projection of $X$ is
  $v_2(\P^{n-k})\subset\P^{N(n-k)}$, then $X$ is
  $v_2(\P^n)\subset\P^{N(n)}$.
\item[{\rm (b)}] If the general $k$-tangential projection of $X$ is a
  projection of $v_2(\P^{n-k})$ into $\P^{N(n-k)-1}$, then $X$ is
  either $v_2(\P^{n})\subset\P^{N(n)-1}$ or $B^n\subset\P^{N(n)-1}$.
\end{enumerate}
\end{theorem}

In order to prove Theorem \ref{thm:main} we use tangential projections
to reduce the problem to smaller dimension, so that an inductive
procedure on the dimension of $X \subset \P^N$ can be achieved.
Another ingredient in this proof is Zak's bound on $N$ in terms of the
secant defect of $X \subset \P^N$ (cf. \cite[Ch.  VI]{Zak}). Note that
we reobtain this bound in Section \ref{sec:drop} as a consequence of
the basic properties of tangential projections.

Let us also remark that our proof of Theorem \ref{thm:USA} does not
rely on Theorem \ref{thm:zak-severi} but we reprove it in a different
way. In particular we have avoided the use of the smoothness of the
so-called {\it entry loci} of $X \subset \P^N$, which is necessary in
Zak's proof of Theorem \ref{thm:zak-severi}.

The structure of the paper is the following: We start by recalling the
notions of secant defects and tangential projections of a projective
variety in Section \ref{sec:prelim}. In Section \ref{sec:drop} we
introduce the {\it drop sequence} of a projective variety $X$, i.e.
the sequence of coranks of the successive tangential projections of
$X$. In Section \ref{sec:main} we develop the proof of Theorem
\ref{thm:main} using the tools described in Sections \ref{sec:prelim}
and \ref{sec:drop}. Finally, in Section \ref{sec:USA}, we obtain
Theorem \ref{thm:USA} as a consequence of Theorems \ref{thm:main} and
\ref{thm:severi-scorza}.


\section{Preliminaries}\label{sec:prelim}
We begin this section by recalling the definition of secant defects
of a projective variety. Subsection \ref{ssec:tangproj} deals with
the definition and basic properties of tangential projections.

\subsection{Secant varieties and defects}\label{ssec:secant}
Throughout the paper $X\subset\P^N$ will denote a complex irreducible
projective variety of dimension $n$.

We consider the sequence of secant varieties of $X$, that is,
\[
X\subsetneq S^1X=SX \subsetneq S^2X\subsetneq \dots\subsetneq S^{k_0}X
\]
where the {\it $k$-secant variety} is defined as:
\[
S^k X=\overline{\{z \in \langle x_0,\dots,x_k \rangle\mid
  (x_0,\dots,x_k)\in U \}},
\]
being $\langle x_0,\dots,x_k\rangle\subset\P^N$ the linear span of the
points $x_0,\dots,x_k\in X$, $U=\{(x_0, \dots, x_k) \in
X^{k+1}\mid\dim(\langle x_0,\dots,x_k \rangle)=k\}$ and $k_0$ the
least integer such that $S^{k_0}X=\langle X \rangle\subset\P^N$. The
expected dimension of $S^k X$ is $(k+1)n+k$ and for $k \leq k_0$ we
denote the difference with the actual dimension by $\delta_k(X)$ and
we call it the {\it $k$-secant defect of $X \subset \P^N$}. We also
set $\delta_k(X)=0$ for $k\leq 0$. If $\delta_k(X)>0$ then $X \subset
\P^N$ is said to be {\it $k$-defective}. We will write $\delta_k$
instead of $\delta_k(X)$ when there is no ambiguity.

As said in the introduction, the study of $1$-defective, not
necessarily smooth, varieties of small dimension goes back to Severi
\cite{severi} and Scorza \cite{scorza} (see also \cite{fujita} and
\cite{ChC}), who completed the classification for dimension two and
three, respectively.

\begin{theorem}\label{thm:severi-scorza} Let $X\subset\P^N$ be a 
  non-degenerate $1$-defective projective variety.

\begin{enumerate}
\item[{\rm (a)}] If $\dim(X)=2$ and $N \geq 5$, then $X$ is either a
  cone over a curve or the Veronese surface $v_2(\P^2)\subset \P^5$.
  
\item[{\rm (b)}] If $\dim(X)=3$ and $N \geq 7$, then one of the
  following holds:
\begin{itemize}
\item[{\rm (i)}] $X$ is a cone,
\item[{\rm (ii)}] $X$ lies in a 4-dimensional cone over a curve,
\item[{\rm (iii)}] $X\subset\P^7$ is contained in a $4$-dimensional
  cone over the Veronese surface $v_2(\P^2)\subset\P^5$,
\item[{\rm (iv)}] $X=v_2(\P^3)\subset\P^9$ or one of its projections
  into $\P^8$ or $\P^7$,
\item[{\rm (v)}] $X\subset\P^7$ is a hyperplane section of the Segre
  embedding $\P^2\times\P^2\subset\P^8$.
\end{itemize}
\end{enumerate}
\end{theorem}

For $Y \subset X \subset \P^N$ define the {\it relative secant variety
  of $X$ with respect to $Y$} as
\[
S(Y,X)=\overline{\{z \in \langle y,x \rangle\mid y\in Y,\; x \in X, \;
  y\neq x\}}.
\]

We denote by $T_xX\subset\P^N$ the projective tangent space to $X$ at
a point $x\in X$. If $Y$ is contained in the smooth part of $X$ then
the {\it relative tangent variety of $X$ with respect to $Y$} is
defined as
\[
T(Y,X)=\bigcup_{y \in Y}T_yX.
\]

Let us recall the following useful consequence \cite[Ch.~I,
Thm.~1.4]{Zak} of Fulton-Hansen's Theorem \cite{FH}.

\begin{lemma}\label{lem:fultonhansen}
  Let $X \subset \P^N$ be a projective variety and let $Y \subset X$
  be an irreducible closed subset contained in the smooth part of $X$.
  Then either:
\begin{enumerate}
\item[{\rm (a)}] $\dim(T(Y,X))=\dim(Y)+\dim(X)=\dim (S(Y,X))-1$, or
\item[{\rm (b)}] $T(Y,X)=S(Y,X)$.
\end{enumerate}
\end{lemma}

The last definition of this section is the following: for a general $u
\in S^kX$ the {\it entry locus of $u$} is defined as
\[
E_u(X)=\overline{\{x \in X\mid \mbox{ there exists } x' \in S^{k-1}X,
  \; u\in \langle x,x'\rangle\}}.
\]
\begin{remark}\label{rem:entry}
  {\rm Note that a simple count of dimensions shows that
    $\dim(E_u(X))=\delta_k-\delta_{k-1}$ for general $u\in S^kX$.}
\end{remark}

Throughout the paper, given a rational map $\pi:X\to Y$ and closed
sets $C\subset X$, $D\subset Y$, we will denote by $\pi(C)$ and by
$\pi^{-1}(D)$ the strict transforms of $C$ in $Y$ and of $D$ in $X$,
respectively.

\subsection{Tangential projections}\label{ssec:tangproj}

Let us recall the definition of tangential projection. We refer the
interested reader to \cite{russo} for a more detailed account. Let us
remark that tangential projections have been used in other problems
regarding projective varieties with special properties on their
projections (see, for instance, \cite{bronowski}, \cite{ChC},
\cite{Ch-Ci}, \cite{CMR}, \cite{chiantinisurvey}, \cite{Ci-Ru}).

\begin{definition}\label{defin:tangproj} Consider
  the notations of Section \ref{ssec:secant}. Given $k \leq k_0$ and
  $(x_1,\dots, x_k) \in U$ general, $\pi_k:X\to X_k$ stands for the
  linear projection of $X$ onto its image $X_k$ from the linear space
  $\langle T_{x_1}X,\dots,T_{x_k}X \rangle$, and we call it the {\it
    $k$-tangential projection of} $X \subset \P^N$. A $1$-tangential
  projection is simply called {\it tangential projection}.
\end{definition}

The following lemma shows how tangential projections can be applied to
compute the dimension of the secant varieties (cf.  \cite{terracini}):

\begin{lemma}[Terracini's Lemma]\label{lem:terracini}
  Let $X \subset \P^N$ be a projective variety and $u\in S^kX$ a
  general point in a general $(k+1)$-secant $k$-space $\langle
  x_0,\dots,x_k\rangle$. Then
\[
T_u S^kX=\langle T_{x_0}X,\dots,T_{x_k}X\rangle.
\]
\end{lemma}

In particular, for $k\leq k_0$ it holds that $\dim(\sga T_{x_0}X,\dots
,T_{x_k}X\sgc)=(k+1)n+k-\delta_k$. This equality has a counterpart in
the relative position of tangent spaces to $X$. If, for example,
$\delta_1=1$ then the tangent spaces to $X$ at two general points meet
in just one point.

The following lemma is a direct consequence of Lemma
\ref{lem:terracini}:

\begin{lemma}\label{dimfiber} Let $X \subset \P^N$ be a projective 
  variety, $k \leq k_0$, and $F^k$ be the general fiber of the
  $k$-tangential projection $\pi_k:X\to X_k$. Then
  $\dim(F^k)=\delta_k-\delta_{k-1}$.
\end{lemma}

\begin{proof}
  By Lemma \ref{lem:terracini}, $\langle T_{x_1}X,\dots,
  T_{x_{k}}X\rangle=T_uS^{k-1}X$ for general $u \in \langle
  x_1,\dots,x_{k}\rangle$. If $x\in X$ is general, we have
  $\dim(\langle T_xX, T_uS^{k-1}X \rangle)=(k+1)n+k-\delta_{k}$. Then
  $\dim(X_k)=(k+1)n+k-\delta_{k}-(kn+k-1-\delta_{k-1})-1=
  n-(\delta_k-\delta_{k-1})$ and so $\dim(F^k)=\delta_k-\delta_{k-1}$.
\end{proof}

The following lemma studies when a general tangent space to $X$
intersects $X$ in codimension $1$. This will be used in the proof of
Theorem \ref{thm:main}. The classical reference to this result is
\cite{DP}. See also \cite[Prop.~5.2]{CMR} where smoothness of $X$ (at
least in codimension two) is required in the proof.

\begin{lemma}\label{lem:basecomponents} Let $X \subset \P^N$ be a 
  non-linear projective variety of dimension $n\geq 2$. For a smooth
  point $x\in X$ define $D_x$ as the $(n-1)$-dimensional part of the
  scheme $X \cap T_xX$.  If $D_x\neq\emptyset$ for general $x \in X$,
  then $X\subset\langle X\rangle$ is either a hypersurface or swept
  out by linear spaces of dimension $n-1$.
\end{lemma}

\begin{proof}
  Cutting with a general $\P^{n-2}$ we reduce the statement to the
  case of surfaces. Consider the family ${\cal F}=\{D_x|\; x\in U\}$
  where $U$ is the open subset of the smooth part of $X$ where $D_x$
  is non-empty (in fact $U$ is the smooth part of $X$ by
  semicontinuity). Let us observe that $\dim({\cal F})>0$, being $X$
  non-linear.
  
  If $\dim({\cal F})=1$ then for $x \in U$ there exists a curve
  $L_x\subset U$ such that $D_x=D_z$ for any $z \in L_x$, whence
  $\bigcap_{z \in L_x} T_zX\supset D_x$. Hence either $D_x$ is a line
  or $T_zX=T_xX$ for any $z \in L_x$. In the first case, the lines
  parameterized by ${\cal F}$ sweep out $X$, otherwise the general
  $T_xX$ contains a fixed line, contradicting the non-linearity of
  $X$. If the latter holds then the reduced structure of $D_x$ is
  linear as a consequence of the linearity of the general fiber of the
  Gauss map (see, for instance, \cite[Ch.~I,~Thm.~2.3]{Zak}).
  
  If $\dim({\cal F})=2$ then for general $z,z' \in X$ there exists
  $x\in X$ such that $z,z' \in D_x$. Since $D_x$ is a plane curve then
  either $D_x=\P^1$ (and so $X=\P^2$, a contradiction), or $D_x$ is a
  plane conic, or $\deg (D_x)>2$. If the latter holds then the general
  secant line to $X$ is trisecant so that $X \subset \P^3$ by the well
  known Trisecant Lemma (cf. \cite[p. 110]{acgh}). If $D_x$ is a plane
  conic then $X$ is a projective surface with a two dimensional family
  of plane conics so that either $X \subset \P^3$, or $X=v_2(\P^2)
  \subset\P^5$, or one of its projections into $\P^4$ (cf.
  \cite{segrecurvasplanas}). But these cases can be excluded because
  $D_x$ is not a conic for general $x\in X$.
\end{proof}

\section{The drop sequence and the defective sequence of a 
  projective variety}\label{sec:drop}

Let $k\leq k_0$ be a positive integer. The general $k$-tangential
projection can be written as a composition of $1$-tangential
projections in the following way. Given a sequence of general points
$x_1,\dots,x_k\in X$ we consider the corresponding sequence of
tangential projections:
\[
\xymatrix{ X \ar[r]_{p_1} \ar@/^0.5cm/[rrr]^{\pi_k}& X_1\ar[r]_{p_2}&
  \dots\ar[r]_{p_k}&X_k},
\]
where $p_1=\pi_1$ is the tangential projection from $T_{x_1}X$ and
$p_{j+1}$ denotes the tangential projection of $X_j$ from
$T_{\pi_j(x_{j+1})}X_j$. Observe that $\pi_j=p_j\circ\dots \circ p_1$,
for all $1 \leq j \leq k_0$, and $X_{k_0}$ is linear.

\begin{definition}\label{def:drop}
  Let $X\subset\P^N$ be a projective variety. The sequence
  $\zeta(X)=\zeta:=(\zeta_1,\dots,\zeta_{k_0})$, where $\zeta_j$ is
  the corank of $p_j$, is called the {\it drop sequence} of $X$. Note
  that $\zeta_1=\delta_1$. The sequence
  $\delta(X)=\delta:=(\delta_1,\dots,\delta_{k_0})$ is called {\it
    defective sequence} of $X$.
\end{definition}

Let us introduce the following notation. Given $a=(a_1,\dots,a_r)$ a
sequence of integers we will write
$da:=(a_1-0,a_2-a_1\dots,a_r-a_{r-1})$ for the sequence of first
differences of $a$.

\begin{remark}\label{rem:2dif}
  {\rm By Lemma \ref{dimfiber} the relation between the drop sequence
    and the defective sequence is $\zeta=d^2(\delta)$.}
\end{remark}

\begin{example}\label{ex:dropveronese}
{\rm  A direct computation for $v_2(\P^n)\subset \P^{N(n)}$ shows that
  $k_0=n$ and $\zeta=(1,\dots,1)$ and for its projection from a point
  we get $k_0=n-1$ and $\zeta=(1,\dots,1)$. See Example
  \ref{ex:scorzageneralized} for further examples.}
\end{example}

\begin{remark}
  {\rm For any sequence of non-negative integers $z=(z_1,\dots,z_r)$
    there exists a projective variety $X\subset\P^N$ such that
    $\zeta=z$ and $k_0=r$, as proved in \cite{CJ}.}
\end{remark}

In the following subsection we recall some arithmetic properties of
the defective sequence of a smooth projective variety.

\subsection{Additivity and superadditivity of the defective 
  sequence of smooth projective varieties}\label{ssec:addit}

We begin by using Lemma \ref{lem:fultonhansen} to prove that, in the
smooth case, $\delta_1$ cannot decrease by linear projections.

\begin{proposition}\label{prop:projections}
  Let $X\subset\P^N$ be a smooth projective variety, $V\subset\p^N$ a
  linear subspace, $p_V$ the linear projection from $\P^N$ with center
  $V$, and $\pi_V:X\to Z$ the corresponding rational map onto its
  image. If $Z$ is not linear, then
\[
\delta_1(Z)\geq\delta_1(X).
\]
\end{proposition}

\begin{proof}
  For general $z\in Z$ we set $X_z:=\pi_V^{-1}(z)$. We first claim
  that $T(X_z,X)\subsetneq S(X_z,X)$. In fact, $T(X_z,X)$ is contained
  in the linear space $p_V^{-1}(T_zZ)$. On the other side, if
  $X\subset p_V^{-1}(T_zZ)$, then $\pi_V(X)=Z\subset T_zZ$
  contradicting the non-linearity of $Z$. It follows from Lemma
  \ref{lem:fultonhansen} applied to each irreducible component of
  $X_z$ that $S(X_z,X)$ has the expected dimension
  $2\dim(X)-\dim(Z)+1$.
  
  Considering the incidence variety:
\[
\xymatrix{ I:=\big\{(z,u)\mid u\in S(X_z,X)\big\}\subset Z\times SX
  \ar[r]^{\hspace{2.5cm}p_2} \ar[d]^{p_1}&SX,\\Z& }
\]
a dimension count tells us that $\dim (p_1(p_2^{-1}(u)))=\delta_1(X)$
for a general $u\in SX$. It follows that $p(u)\in SZ$ and its entry
locus $E_{p(u)}(Z)$ contains $p_1(p_2^{-1}(u))$. Therefore
$\delta_1(Z)\geq\delta_1(X)$.
\end{proof}

\begin{remark}{\rm Smoothness cannot be dropped
    in Proposition \ref{prop:projections}. Consider, for instance, a
    $2$-dimensional cone $X\subset\P^N$, $N\geq 6$, whose vertex is a
    point and let $V=T_xX$ for a general $x\in X$. Then
    $\delta_1(X)=1$ and $Z=X_1\subset\P^{N-3}$ is not a plane curve,
    so $\delta_1(Z)=0$.}
\end{remark}

An immediate corollary of this result is what we call {\it
  superadditivity} of the defective sequence (cf. \cite[Ch.~V,
Thm.~1.8]{Zak}, having in mind that definitions of $\delta_k$ do not
coincide; see also \cite{zak2} and \cite{fantechi} for a more general
statement):

\begin{corollary}\label{cor:subaddit}
  Let $X \subset \P^N$ be a smooth projective variety with drop
  sequence $(\zeta_1,\dots,\zeta_{k_0})$. Then $\zeta_i\geq\delta_1$
  for all $i$, and the defective sequence of $X$ verifies the {\rm
    superadditivity property} $\delta_k\geq\delta_{k-1}+k\delta_1$. In
  particular if $\delta_1>0$, then $k_0\leq\frac{n}{\delta_1}$.
\end{corollary}

\begin{proof}
  Note that $\pi_k(X)$ is linear if and only if $k=k_0$. By
  Proposition \ref{prop:projections} applied to the linear projection
  $\pi_k$, $\zeta_k=\delta_1(X_{k-1})\geq \delta_1(X)$ for $k\leq
  k_0$. We conclude by noting that
  $\delta_k=\delta_{k-1}+\sum_{i=1}^k\zeta_i$ (see Remark
  \ref{rem:2dif}).
\end{proof}

In this note we are interested in varieties whose defective sequence
satisfy a stronger condition, that we call \emph{additivity}.

\begin{definition}\label{def:addit}
  {\rm We say that the defective sequence of a projective variety
    $X\subset\P^N$ is {\it additive} when
    $\delta_k=\delta_{k-1}+k\delta_1$ (equivalently
    $\delta_k=\delta_1\frac{k(k+1)}{2}$) for every $k\in
    \{1,\dots,k_0\}$, or, in other words, when the drop sequence
    $\zeta$ of $X$ is constant (equivalently $\dim (F^k)=k\delta_1$).}
\end{definition}

\begin{example}\label{ex:scorzageneralized}
  {\rm The defective sequence $\delta$ is additive for the following
    $1$-defective varieties:
\begin{enumerate}
\item[(i)] The Veronese embedding $X=v_2(\P^n)\subset\P^{N(n)}$. Note
  that $X_k=v_2(\P^{n-k})\subset\P^{N(n-k)}$, $\delta_1=1$ and
  $k_0=n$.
  
\item[(ii)] The projected Veronese embedding
  $X=v_2(\P^n)\subset\P^{N(n)-1}$ (resp. $X=B^n\subset\P^{N(n)-1}$).
  Now $X_k=v_2(\P^{n-k})\subset\P^{N(n-k)}$ (resp.
  $X_k=B^{n-k}\subset\P^{N(n-k)}$), $\delta_1=1$ and $k_0=n-1$.
  
\item[(iii)] The Segre embedding
  $X=\P^a\times\P^b\subset\P^{(a+1)(b+1)-1}$. Here
  $X_k=\P^{a-k}\times\P^{b-k}\subset\P^{(a+1-k)(b+1-k)-1}$,
  $\delta_1=2$ and $k_0=\min \{a,b\}$.
  
\item[(iv)] The Pl\"ucker embedding $X=G(1,r)\subset\P^{N(r-1)}$ of
  the grassmannian of lines in $\P^r$. In this case $X_k=G(1,r-2k)$,
  $\delta_1=4$ and $k_0=\frac{r}{2}-1$ if $r$ is even or
  $k_0=\frac{r-1}{2}$ if $r$ is odd.
  
\item[(v)] The Cartan variety $X=E^{16}\subset\P^{26}$, where
  $\delta_1=8$ and $k_0=2$.
\end{enumerate}}
\end{example}

\begin{remark}\label{rem:open}
  {\rm Examples (i), (iii) with $|a-b|\leq 1$, (iv) and (v) are the
    so-called \emph{Scorza varieties} (see \cite[Ch.~VI]{Zak}). This
    list of examples shows that Scorza varieties are contained in the
    class of varieties verifying the more general property of
    additivity on its defective sequence. This suggests that a natural
    further development of this theory is the classification of these
    varieties.}
\end{remark}

\begin{remark}
  {\rm In examples (i)-(iv) we get non-finite sequences of varieties
    $\{X^j\}_{j\in\N}$ such that for any $j \in \N$ the $k$-tangential
    projection of $X^j$ verifies $X^j_{k}=X^{j-k}$. It would be of
    interest to find some other examples of these sequences.}
\end{remark}

If the codimension of $X\subset\P^N$ is big enough, the defective
sequence of $X$ verifies {\it additivity}.

\begin{lemma}\label{lem:addit}
  Let $X\subset\p^N$ be a non-degenerate smooth projective variety of
  dimension $n$. Then
\begin{equation}\label{eq:addit}
N\leq\phi(n,k_0,\delta_1):=n(k_0+1)-k_0
(\delta_1-1)-\delta_1\frac{k_0(k_0-1)}{2}.
\end{equation}

Moreover, equality holds if and only if $\delta$ is additive.
\end{lemma}

\begin{proof}
  Let $d$ denote the dimension of $X_{k_0}$. It follows from Lemma
  \ref{dimfiber} and Corollary \ref{cor:subaddit} that
\begin{equation}\label{eq:dim1}
d=n-\dim (F^{k_0})=n-(\delta_{k_0}-\delta_{k_0-1})\leq
n-k_0\delta_1.
\end{equation}
Since $X\subset\P^N$ is non-degenerate $X_{k_0}$ coincides with its
linear span, whence
\begin{equation}\label{eq:dim2}\begin{array}{rl}
\vspace{0.2cm}d=& N-\dim(S^{k_0-1}X)-1=N-(k_0n+k_0-1-
\delta_{k_0-1})-1\\\geq &
N-(n+1)k_0+\delta_1\frac{k_0(k_0-1)}{2}.\end{array}
\end{equation}
Joining (\ref{eq:dim1}) and (\ref{eq:dim2}) we get the desired result.
For the second assertion, note that equality in (\ref{eq:dim1}) and
(\ref{eq:dim2}) holds if and only if $\delta$ is additive.
\end{proof}

As a by-product of the previous lemma we obtain the following well
known bound (cf. \cite[Ch.~V, Thm.~2.3]{Zak}).

\begin{corollary}\label{lem:bound}
  Let $X \subset \P^N$ be a non-degenerate smooth projective variety
  of dimension $n$. Assume $\delta_1>0$, and let $r_0$ be the rest of
  $n$ modulo $\delta_1$. Then:
\[
N\leq
\frac{1}{2\delta_1}\big(n(n+\delta_1+2)+r_0(\delta_1-r_0-2)\big).
\]

In particular, $N \leq N(n)$ and this bound is sharp.
\end{corollary}

\begin{proof}
  Note that, once $n$ and $\delta_1$ are fixed, the maximum of
  $\phi(n,k_0,\delta_1)$ is achieved at $k_0=(n-r_0)/\delta_1$. A
  simple computation provides the claimed upper bound. For the
  sharpness just consider $v_2(\P^n) \subset \P^{N(n)}$.
\end{proof}

If we assume additivity of the defective sequence, we get some
restrictions on the singularities of the $k$-tangential projections of
$X$:

\begin{lemma}\label{lem:conesproj}
  Let $X \subset \P^N$ be a smooth projective variety. If $\delta$ is
  additive, then $X_k$ and $S^kX$ are not cones for $k<k_0$.
\end{lemma}

\begin{proof}
  Assume $X_k$ is a cone with vertex $V_k$ over a variety $X'_k$. The
  hypothesis $k<k_0$ implies that $X_k$ is not linear, whence $X'_k$
  is not linear. Moreover $X'_k$ is the linear projection of $X_k$
  from $V_k$, so it is also a linear projection of $X$. Then
  Proposition \ref{prop:projections} implies that
  $\delta_1(X'_k)\geq\delta_1(X)$. But using that
\[
\delta_1(X_k)=\delta_1(X'_k)+\dim (V_k)+1,
\]
and that $\zeta_1=\zeta_{k+1}$ (whence $\delta_1(X)=\delta_1(X_k)$)
since $\delta(X)$ is additive, we get the contradiction
\[
\delta_1(X)=\delta_1(X_k)>\delta_1(X'_k)\geq\delta_1(X).
\]
For $S^kX$ not to be a cone we have to prove that the following set is
empty:
\[
V:=\bigcap_{\scriptsize{\begin{array}{c}x_1,\dots,x_{k+1}\in
      X\\\mbox{general}\end{array}}}\langle
T_{x_1}X,\dots,T_{x_{k+1}}X\rangle.
\]
But since $X_r$ is not a cone, then
\[
\bigcap_{x_{r}\in X}\langle T_{x_1}X,\dots,T_{x_{r}}X\rangle=\langle
T_{x_1}X,\dots,T_{x_{r-1}}X\rangle,\mbox{ for all
}r\in\{2,\dots,k+1\}.
\]
Recursively we obtain $V=\bigcap_{x_1\in X}T_{x_1}X$, which is empty
since $X$ itself is not a cone.
\end{proof}

\begin{remark}
  {\rm If $X \subset \P^N$ is a smooth variety but $\delta$ is not
    additive, then $X_k$ might be a cone for $k<k_0$. Consider, for
    instance, an integer $q\geq 4$ and a rational normal scroll
    $S_{1,q}\subset\P^{q+2}$. Then $X_1=S_{0,q-2}\subset\P^{q-1}$ is a
    non-linear cone.}
\end{remark}

\section{Proof of Theorem \ref{thm:main}}\label{sec:main}

First, we reduce the proof to the case $k=n-2$.

\begin{lemma}\label{lem:k=n-2}
  If Theorem \ref{thm:main} holds for $k=n-2$, then it holds for any
  $k<k_0$.
\end{lemma}

\begin{proof}
  Recall that $k_0\leq n$ by Corollary \ref{cor:subaddit}.
  
  Consider first $k<n-2$. Assume that $X_k$ is either
  $v_2(\P^{n-k})\subset\P^{N(n-k)}$ or one of its projections into
  $\P^{N(n-k)-1}$. Then $X_{n-2}=\pi_{n-k-2}(X_k)$ is either
  $v_2(\P^{2})\subset\P^5$ or one of its projections into $\P^4$,
  respectively. Since by hypothesis Theorem \ref{thm:main} holds for
  $k=n-2$, then $X$ is either $v_2(\P^{n})\subset\P^{N(n)}$, or
  $v_2(\P^{n})\subset\P^{N(n)-1}$, or $B^{n}\subset\P^{N(n)-1}$.
  
  Finally, if $k=n-1$ then necessarily $X_{n-1}=v_2(\P^1)\subset\P^2$,
  so $X\subset\P^{N(n)}$ and its secant defect is additive by Lemma
  \ref{lem:addit}. It follows from Proposition \ref{prop:projections}
  that $X_{n-2}\subset\P^5$ is not a cone, whence
  $X_{n-2}=v_2(\P^2)\subset\P^5$ by Theorem \ref{thm:severi-scorza}(a)
  and so $X=v_2(\P^{n})\subset\P^{N(n)}$.
\end{proof}

The following result is the heart of the paper.

\begin{theorem}\label{prop:main}
  Theorem \ref{thm:main} holds for $k=n-2$.
\end{theorem}

\begin{proof}
  We present the proof divided in several steps.
  
  \medskip

\noindent \emph{Set up}: For general points $x_1,\dots,x_{n-2}\in X$
denote $T:=\langle T_{x_1}X,\dots, T_{x_{n-2}}X \rangle$, so that
$\pi_T:X\to X_{n-2}$ is the projection of $X$ from $T$ onto $X_{n-2}$.
By hypothesis there exists a birational morphism
$\alpha:X_{n-2}\to\p^2$. Denote by $M$ the linear system
$\alpha^*|{\cal O}_{\P^2}(1)|$ on $X_{n-2}$ and observe that the
general $Q \in M$ is a conic $Q\subset X_{n-2}$. Let $L(T)$ be the
$2$-dimensional base component free linear system on $X$ defining the
rational map $\alpha\circ\pi_T$. Finally $F^T_Q\subset X$ stands for
the element of $L(T)$ corresponding to $Q$.

\medskip

\noindent \emph{Step 1: $x_1, \dots,x_{n-2} \in F^T_Q$ for every
  $Q\in M$.}

It suffices to show that $T \cap X$ has no $(n-1)$-dimensional
components meeting $\{x_1,\dots, x_{n-2}\}$. Suppose on the contrary
that there exists an $(n-1)$-dimensional irreducible component
$D\subset X\cap T$ through one of the points, say $x_1$.

Denote $T':=\langle T_{x_{2}}X,\dots,T_{x_{{n-2}}}X\rangle$ and
$\pi_{T'}:X \to X_{n-3}$ the corresponding projection. Since
$x_1,\dots,x_{n-2}$ are taken general and $X\nsubseteq T$, then
$T_{\pi_{T'}(x_{1})}X_{n-3}\cap X_{n-3}$ has a $2$-dimensional
component $\pi_{T'}(D)$ through $\pi_{T'}(x_{1})$. As $X_{n-3}$ is not
a hypersurface (otherwise $k_0=n-2=k$), Lemma \ref{lem:basecomponents}
implies that it is swept out by planes.  Hence $X_{n-2}$ is swept out
by lines, and so $X_{n-2}=B^2\subset\P^4$. Then, by Theorem
\ref{thm:severi-scorza}(b), $X_{n-3}\subset\P^8$ is either $B_3$,
contradicting Lemma \ref{lem:basecomponents}, or a cone, contradicting
Lemma \ref{lem:conesproj}. This concludes Step 1.

\medskip

\noindent \emph{Step 2: $F:=F_Q^T$ is $1$-defective for general
  $Q\in M$, and $\delta_1(F_Q^T)=1$.}

First we prove that $\delta_1(X)=1$. By Lemmas \ref{dimfiber} and
\ref{lem:addit} we get
\[
n-2=\delta_{n-2}(X)-\delta_{n-3}(X)\geq (n-2)\delta_1(X),
\]
so $\delta_1(X)=1$ and for $y_1,y_2 \in F$ general points
$\dim(\langle T_{y_1}X,T_{y_2}X\rangle)=2n$ and $\dim (T_{y_1}X\cap
T_{y_2}X)=0$.

Now we claim that $T_{y_1}F\cap T_{y_2}F = T_{y_1}X\cap T_{y_2}X$.
Denote $r=\dim(\langle T_{y_1}F,T_{y_2}F\rangle)$ and
$s=\dim(T\cap\langle T_{y_1}F,T_{y_2}F\rangle)$. Clearly $T_{y_1}F\cap
T_{y_2}F\subset T_{y_1}X\cap T_{y_2}X$, whence $\dim (T_{y_1}F\cap
T_{y_2}F)\leq 0$ or, equivalently, $r\geq 2n-2$.

On the other hand $\dim(\langle
T_{\pi_T(y_1)}X_{n-2},T_{\pi_T(y_2)}X_{n-2}\rangle)=4$, so we get
$\dim(T\cap\langle T_{y_1}X,T_{y_2}X\rangle)=2n-5.$ Since
$\dim(\langle T_{\pi_T(y_1)}Q,T_{\pi_T(y_2)}Q\rangle)=2$ it follows
that $r-s-1=2$. Finally $T\cap\langle T_{y_1}F,T_{y_2}F\rangle\subset
T\cap\langle T_{y_1}X,T_{y_2}X\rangle,$ then $s\leq 2n-5$ and so
$r\leq 2n-2$.  Therefore $r=2n-2$ and $\dim (T_{y_1}F\cap
T_{y_2}F)=0$, proving the claim.

\medskip

\noindent \emph{Step 3: The complete linear system $L'(T)$
  containing $L(T)$ does not depend on $T$.}

We prove it first in the case $n=3$. By assumption the general
tangential projections of $X$ are isomorphic and we will identify
them, via a fixed isomorphism, with a given one $X_1$. Let $x$ and $y$
be two general points in $X$ and we want to prove that
$L'(T_xX)=L'(T_yX)$. Note that $F_Q^{T_xX}$ and $F_{Q}^{T_yX}$ are
algebraically equivalent for any $Q\in M$. Hence it is enough to prove
that there exists $F\in L'(T_xX)$ such that ${\pi_{T_yX}(F)}\in M$. We
choose $Q_x\in M$ containing $\pi_{T_xX}(y)$ and take
$F:=F_{Q_x}^{T_xX}$. Note that $x\in F$ by Step 1.  Since $y$ is
general, Step 2 implies that $\dim(T_yF\cap T_zF)=0$ for general $z\in
F$, therefore $\dim(\pi_{T_yX}(F))=n-2$.  Moving $y\in F$ we construct
an algebraic family of divisors in $X_1$ containing
$\pi_{T_xX}(F)=Q_x$ and $Q_y:=\pi_{T_yX}(F)$. Since algebraic and
linear equivalence of divisors in $X_1$ coincide, it follows that
$Q_y\in M$.

Now we use recursively the above argument in order to prove the
general case. Let us write
\[
T:=\langle T_{x_1}X,\dots,T_{x_{n-2}}X\rangle,\quad T':=\langle
T_{y_1}X,\dots,T_{y_{n-2}}X\rangle,
\]
where $x_1, \dots, x_{n-2}, y_1, \dots y_{n-2}$ are general points of
$X$. Denote by $X_{n-3}$ the $(n-3)$-tangential projection of $X$ from
$\langle T_{x_1}X,\dots,T_{x_{n-3}}X\rangle$ and consider the
tangential projections of $X_{n-3}$ corresponding to $x_{n-2}$ and
$y_{n-2}$. Note that this two points are smooth in $X_{n-3}$. Since
$X_{n-3}$ is not a hypersurface and it is not swept out by planes, we
conclude that $x_{n-2}\in\Bs(L'(T_{x_{n-2}}X_{n-2}))$. Thus we argue
as in the previous paragraph to deduce that
\[
L'(T)=L'(\langle T_{x_1}X,\dots,T_{x_{n-3}}X,T_{y_{n-2}}X\rangle).
\]
Recursively we obtain the desired result:
\[
\begin{array}{l}
L'(T)=L'(\langle
  T_{x_1}X,\dots,T_{x_{n-3}}X,T_{y_{n-2}}X\rangle)=L'(\langle
  T_{x_1}X,\dots,T_{y_{n-3}}X,T_{y_{n-2}}X\rangle)=\\=\dots=L'(\langle
  T_{y_1}X,\dots,T_{y_{n-3}}X,T_{y_{n-2}}X\rangle)=L'(T').
\end{array}
\]
Along the rest of the proof we will write $L':=L'(T)$ and
\[
L:=\overline{\bigcup_{T} L(T)}\subset L'.
\]

\noindent \emph{Step 4: $L$ is a base point free linear system and
  the general $F_Q^T$ is smooth and irreducible.}

First of all, we claim that $\dim (L)=n$.  Consider the incidence
variety:
\[
\xymatrix{I=\left\{(x_1,\dots,x_{n-2},F)\mid F\in L(T)\right\}\subset
  V\times L\ar[r]^{\hspace{3.1cm} p_2}
  \ar[d]^{p_1} & L\\
  V}
\]
where $V\subset X^{n-2}$ is the non-empty open subset defined by
$(n-2)$-uples in the hypotheses of the Theorem. We have $\dim(
V)=(n-2)n$ and $\dim (p_1^{-1}(v))=2$ for each $v\in V$. This implies
that $\dim (I)=(n-2)n+2$. Furthermore $\dim(p_2^{-1}(F))=(n-2)(n-1)$
for general $F\in L$, since we have shown in the previous step that
the image of $F$ by the $(n-2)$-tangential projection of $X$ is a
conic for a general choice of points $x_1,\dots,x_{n-2}\in F$. Hence
$\dim( L)=(n-2)n+2-(n-2)(n-1)=n$.

Now we prove that $L\subset L'$ is linear. By \cite{bsegre} it
suffices to show that $L$ contains a $3(n-2)$-dimensional family $\cG$
parameterizing the planes $L(T)$. A dimension count shows that
$\dim(\cG)=\dim(V)-(n-2)\dim(X_T),$ where $X_T=\{x \in X\mid T_xX
\subset T\}$. We claim that $\dim(X_T) \leq n-3$. Consider $X_{n-3}$
the $(n-3)$-tangential projection of $X$ from general points
$x_{1},\dots,x_{n-3} \in X_T$ and $F^{n-3}$ its general fiber.
Observe that $\dim (F^{n-3})=\delta_{n-3}-\delta_{n-4}=n-3$ by Lemma
\ref{dimfiber} and Lemma \ref{lem:addit}, having in mind that
necessarily $1=\zeta_1=\dots=\zeta_{n-2}$. If
$\dim(X_T)>n-3=\dim(F^{n-3})$ then $X_{n-3}$ is developable, that is
the general tangent space is tangent along a subvariety of positive
dimension, and not a cone by Lemma \ref{lem:conesproj}. Hence, by
Theorem \ref{thm:severi-scorza}, $X_{n-3}$ is contained in a
$4$-dimensional cone over a curve so that $X_{n-2}$ is a cone,
contradicting again Lemma \ref{lem:conesproj}. Hence $\dim(\cG) \geq
(n-2)n-(n-2)(n-3)=3n-6,$ as claimed, so that $L\subset L'$ is linear.

Once we know that $L$ is linear Step 3 implies that $\Bs (L)
\subset\bigcap_{T}\Bs L(T)\subset\bigcap _{u \in
  S^{k-1}X}T_uS^{k-1}X=\emptyset$ by Lemma \ref{lem:conesproj}. Hence
$L$ is base point free and we get the first assertion in the statement
of Step 4. Finally, applying Bertini's Theorems \cite[II Thm.~8.18 and
III Ex.~11.3]{Hartshorne} to $L$ we get the smoothness and
irreducibility for the general $F\in L$.  \medskip

\noindent \emph{Step 5: End of the proof}

We prove the theorem by induction on $n$. If $n=3$ then the result is
a consequence of Theorem \ref{thm:severi-scorza}(b).

We prove it for $n$. Let $F\in L$ be a general element and
$x_1,\cdots,x_{n-2}$ general points in $F$.  By Steps 4 and 2 $F$ is
smooth, irreducible and $\delta_1(F)=1$. In particular its defective
sequence $\delta(F)$ is superadditive by Corollary \ref{cor:subaddit}.
Denote $T'=\langle T_{x_1}F,\dots,T_{x_{n-2}}F\rangle$. We claim that
$T'=T\cap \langle F\rangle$ and in particular the $(n-2)$-tangential
projection of $F$ from $T'$ coincides with $\pi_T|_F$. Note first that
$T'\subset T\cap\langle F\rangle$, and so $\pi_T|_F$ factors through
$\pi_{T'}$. Since $\dim(\pi_T(F))=1$ it follows that $\delta(F)$ is
additive.  By Lemma \ref{lem:addit}, since $\delta_1(F)=1$ and
$k_0(F)=n-1$ then $\dim (\langle F \rangle) = N(n-1)$. Hence
$F_{n-3}=v_2(\P^2)\subset\P^5$ as in the proof of Lemma
\ref{lem:k=n-2}. By induction $F=v_2(\P^{n-1}) \subset \P^{N(n-1)}$.
The linear system $L$ on $X$ defines a map $\varphi:X \to \P^n$. Note
that $\varphi$ is a birational map since its restriction
${\varphi}_{|F}:F=v_2(\P^{n-1})\to\P^{n-1}$ is an isomorphism. The
birational inverse $\varphi^{-1}:\P^n\to X$ is defined by a linear
subspace $|V|$ of $|\cO_{\P^n}(a)|$. Since
$\varphi^{-1}_{|\P^{n-1}}:\P^{n-1}\to F$ is given by
$|\cO_{\P^{n-1}}(2)|$, it follows that $a=2$.

If $X_{n-2}=v_2(\P^2)\subset\P^5$ then $|V|=|\cO_{\P^n}(2)|$ by a
count of dimensions, whence $X=v_2(\P^n)\subset\P^{N(n)}$.

If $X_{n-2}\subset\P^4$ is a projection of $v_2(\P^2)\subset\P^5$,
then $|V|$ is a codimension $1$ linear subspace of $|\cO_{\P^n}(2)|$.
Therefore $X\subset\P^N$ is a projection of
$v_2(\P^n)\subset\P^{N(n)}$ from a point. Finally, since $X$ is
smooth, this point is either general in $\P^{N(n)}$ or a point of
$v_2(\P^n)$, as stated in Theorem \ref{thm:main}. \end{proof}

\section{Proof of Theorem \ref{thm:USA}}\label{sec:USA}

In this section we prove Theorem \ref{thm:USA}. We first show that the
defective sequence of $X$ in the hypotheses of the theorem is
additive. This allows us to apply Lemma \ref{lem:conesproj} to the
$(n-2)$-tangential projection of $X$ and the result follows from
Theorem \ref{thm:main}.

\begin{lemma}\label{lem:defects}
  Under the hypotheses of Theorem \ref{thm:USA} it follows that the
  defective sequence of $X$ is additive and either $N=N(n)-1$ and
  $k_0=n-1$, or $N=N(n)$ and $k_0=n$.
\end{lemma}

\begin{proof}
  Using Lemma \ref{lem:addit}, under the hypotheses of Theorem
  \ref{thm:USA} we get
\[
N(n)-1=\frac{n^2+3n-2}{2}\leq N\leq n(k_0+1)-\frac{k_0(k_0-1)}{2},
\]
which implies $k_0=n-1$ or $n$. If $k_0=n-1$ the inequalities above
are equalities and additivity holds by Lemma \ref{lem:addit}. If
$k_0=n$, we have two possibilities: either $N=N(n)$ and additivity
holds, or $N=N(n)-1$. But in the second case $\delta$ is not additive
and the drop sequence is $(1,\dots,1,2)$. Therefore, $X_{n-1}$ is a
curve of defect $2$, hence linear, contradicting $k_0=n$.
\end{proof}

\begin{proof}[Proof of Theorem \ref{thm:USA}] 
  From Lemma \ref{lem:defects}, we get the following two
  possibilities:
  
  If $N=N(n)$ then $X_{n-2}\subset\P^5$ is a defective surface. From
  Theorem \ref{thm:severi-scorza}(a) and Lemma \ref{lem:conesproj} it
  follows that $X_{n-2}=v_2(\P^2)\subset \P^5$. Now Theorem
  \ref{thm:main} applies so that $X=v_2(\p^n)\subset\p^{N(n)}$.
  
  If $N=N(n)-1$ then $X_{n-3}\subset\P^8$ is a defective threefold. It
  follows from Theorem \ref{thm:severi-scorza}(b) that
  $X_{n-3}\subset\P^8$ is either a cone over a surface, or is
  contained in a $4$-dimensional cone over a curve, or is a projection
  of $v_2(\P^3)\subset \P^9$. In the first two cases $X_{n-2}$ is a
  cone, so they are again discarded by Lemma \ref{lem:conesproj}. In
  the latter case $X_{n-2}\subset\P^4$ is a projection of
  $v_2(\P^2)\subset \P^5$. Then we conclude again via Theorem
  \ref{thm:main}.
\end{proof}

Finally we reformulate Theorem \ref{thm:USA} in the following way,
that provides further motivations to the suggestion presented in
Remark \ref{rem:open}:

\begin{corollary} Let $X\subset\P^N$ be a non-degenerate smooth 
  $1$-defective projective variety of dimension $n$. Assume $\delta$
  is additive. Then:
\begin{enumerate}
\item[{\rm (a)}] If $k_0=n>1$, then $X$ is
  $v_2(\p^n)\subset\p^{N(n)}$.
\item[{\rm (b)}] If $k_0=n-1>1$, then $X$ is either
  $v_2(\p^n)\subset\p^{N(n)-1}$ or $B^n\subset\p^{N(n)-1}$.
\end{enumerate}
\end{corollary}

\begin{proof} It follows from Corollary \ref{cor:subaddit} and 
  Lemma \ref{lem:addit} that $N=n(k_0+1)-\frac{k_0(k_0-1)}{2}$. Hence
  $N\geq N(n)-1$.
\end{proof}


\hrule \medskip
\par\noindent
Addresses:

\par\noindent{\tt roberto.munoz@urjc.es}
\par\noindent Departamento de Matem\'atica Aplicada, 
Universidad Rey Juan Carlos, 28933 M\'ostoles (Madrid), Spain
\par\noindent {\tt jcsierra@mat.ucm.es}
\par\noindent Departamento de \'Algebra, Facultad de Matem\'aticas, 
Universidad Complutense de Madrid, Ciudad Universitaria, 28040 Madrid,
Spain.
\par\noindent{\tt luis.sola@urjc.es}
\par\noindent Departamento de Matem\'atica Aplicada, 
Universidad Rey Juan Carlos, 28933 M\'ostoles (Madrid), Spain


\begin{thebibliography}{99}
%
%
  
\bibitem[ACGH]{acgh} E. Arbarello, M. Cornalba, P. Griffiths \and J.
  Harris, \emph{Geometry of algebraic curves}. Grudl. der Math. 267
  (Springer-Verlag, Berlin 1985).
  
\bibitem[ASU]{ASU} E. Arrondo, J.C. Sierra \and L. Ugaglia,
  `Classification of $n$-dimensional subvarieties of $G(1,2n)$ that
  can be projected to $G(1,n+1)$', \emph{Bull. London Math. Soc. }37
  (2005) 673--682.
  
\bibitem[B]{bronowski} J. Bronowski, `The sum of powers as canonical
  expressions', \emph{Proc. Cam. Phil. Soc. }29 (1932) 69--82.
  
\bibitem[CJ]{CJ} M.L. Catalano-Johnson, `The possible dimensions of
  the higher secant varieties', \emph{Amer. J. Math. }118 (1996)
  355--361.
  
\bibitem[Ch]{chiantinisurvey} L. Chiantini, `Lectures on the structure
  of projective embeddings'. \emph{Rend. Sem. Mat. Univ. Politec.
    Torino }62 (2004) 335--388.
  
\bibitem[ChC1]{ChC} L. Chiantini \and C. Ciliberto, `Threefolds with
  degenerate secant variety: on a theorem of G. Scorza'.  Geometric
  and combinatorial aspects of commutative algebra (Messina, 1999),
  111--124, \emph{Lecture Notes in Pure and Appl. Math. }217 (Dekker,
  New York, 2001).
  
\bibitem[ChC2]{Ch-Ci} L. Chiantini \and C. Ciliberto, `Weakly
  defective varieties', \emph{Trans. Amer. Math. Soc. }354 (2002)
  151--178.
  
\bibitem[CMR]{CMR} C. Ciliberto, M. Mella \and F. Russo, `Varieties
  with one apparent double point', \emph{J. Algebraic Geom.  }13
  (2004) 475--512.
  
\bibitem[CR]{Ci-Ru} C. Ciliberto \and F. Russo, `Varieties with
  minimal secant degree and linear systems of maximal dimension on
  surfaces', \emph{Adv. Math. }200 (2006) 1--50.
  
\bibitem[DP]{DP} P. Del Pezzo, `Sugli spazi tangenti ad una superficie
  o ad una variet\`a immersa in uno spazio a pi\`u dimensionali',
  \emph{Rend. Accad. Sci. Fis. Mat. Napoli }25 (1886) 176--180.
  
\bibitem[Fa]{fantechi} B. Fantechi, `On the superadditivity of secant
  defects', \emph{Bull. Soc. Math. France }118 (1990) 85--100.
  
\bibitem[Fu]{fujita} T. Fujita, `Projective threefolds with small
  secant varieties', \emph{Sci. Papers College Gen. Ed. Univ.  Tokyo
  }32 (1982) 33--46.
  
\bibitem[FuR]{fu-ro} T. Fujita \and J. Roberts, `Varieties with small
  secant varieties: the extremal case', \emph{Amer. J. Math.  }103
  (1981) 953--976.
  
\bibitem[FH]{FH} W. Fulton \and J. Hansen, `A connectedness theorem
  for projective varieties, with applications to intersections and
  singularities of mappings', \emph{Ann. of Math. }110 (1979)
  159--166.
  
\bibitem[H]{Hartshorne} R. Hartshorne, \emph{Algebraic geometry}.
  Graduate Texts in Mathematics 52 (Springer-Verlag, New
  York-Heidelberg, 1977).
  
\bibitem[Ru]{russo} F. Russo, \emph{Tangents and secants of algebraic
    varieties: notes of a course}. Publica\c c\~oes Matem\'aticas do
  IMPA. 24$\sp {\rm o}$ Col\'oquio Brasileiro de Matem\'atica
  (Instituto de Matem\'atica Pura e Aplicada (IMPA), Rio de Janeiro,
  2003).
  
\bibitem[Sc1]{scorza} G. Scorza, `Determinazione delle variet\`a a tre
  dimensioni di $S_r$, $r \geq 7$, i cui $S_3$ tangenti si tagliano
  due a due', \emph{Rend. Cir. Mat. Palermo }25 (1908) 913--920.
  
\bibitem[Sc2]{scorza2} G. Scorza, `Sulle variet\`a a quattro
  dimensioni di $S_r$ ($r \geq 9$), i cui $S_4$ tangenti si tagliano
  due a due', \emph{Rend. Cir. Mat. Palermo }27 (1909) 148--178.
  
\bibitem[Se1]{bsegre} B. Segre, `Sulle $V_ n$ contenenti pi\`u di
  $\infty^{n-k}$ $S_k$ I, II', \emph{Atti Accad. Naz. Lincei, VIII.
    Ser., Rend., Cl. Sci. Fis. Mat. Nat. (8) }5 (1948) 275--280.
  
\bibitem[Se2]{segrecurvasplanas} C. Segre, `Le superficie degli
  iperspazi con una doppia infinit\`a di curve piane o spaziali',
  \emph{Atti Accad. Sci. Torino Cl. Sci. Fis. Mat. Natur. }56
  (1920-1921) 75--89.
  
\bibitem[S]{severi} F. Severi, `Intorno ai punti doppi impropri di una
  superficie generale dello spazio a quattro dimensioni e ai suoi
  punti tripli apparenti', \emph{Rend. Cir. Mat.  Palermo }15 (1901)
  33--51.
  
\bibitem[T]{terracini} A. Terracini, `Sulle $V_k$ per cui la variet\`a
  degli $S_h$ $(h+1)$-seganti ha dimensione minore dell'ordinario',
  \emph{Rend. Circ. Mat. Palermo }31 (1911) 392--396.
  
\bibitem[Z1]{zak2} F.L. Zak, `Linear systems of hyperplane sections on
  varieties of small codimension', \emph{Funktsional.  Anal. i
    Prilozhen }19 (1985) 1--10.
  
\bibitem[Z2]{Zak} F.L. Zak, \emph{Tangents and secants of algebraic
    varieties}. Translations of Mathematical Monographs, 127 (American
  Mathematical Society, Providence, RI, 1993).

\end{thebibliography}
\end{document}